\documentclass[12pt,  amstex, amscd]{article}
\usepackage{amsmath,amssymb}
\newtheorem{deff}{Deffinition}[section]

\newtheorem{theo}[deff]{Theorem}
\newtheorem{pro}[deff]{Proposition}

\newtheorem{re}{Remark}

\newtheorem{lem}[deff]{Lemma}

\newcommand{\s}{\ensuremath{\star}}

\renewcommand{\v}{\ensuremath{\wedge}}

\newcommand{\A}{\ensuremath{\tilde A      }  }

\newcommand{\B}{\ensuremath{\tilde B  }  }

\newcommand{\g}{\ensuremath{\frak g}}

\newcommand{\ds}{\ensuremath{\displaystyle}}

\newcommand{\h}{\ensuremath{\frak h}}

\newcommand{\gs}{\ensuremath{{\frak g}^{\ast}}}

\newcommand{\Ci}{\ensuremath{C^{\infty}}}

\newcommand{\pp}{\ensuremath{{p\cdots p}}}

\newcommand{\qq}{\ensuremath{{q\cdots q}}}

\newcommand{\pqp}{\ensuremath{{p\cdots q\cdots p}}}

\newcommand{\qpq}{\ensuremath{{q\cdots p\cdots q}}}

\newcommand{\fp}{\ensuremath{{\cal F}_{p}}}

\newcommand{\la}{\ensuremath{{\lambda}}}

\newcommand{\p}{\ensuremath{\partial}}

\begin{document}
\title{Deformation quantization and quantum coadjoint orbits of SL(2,$\mathbb R$).}         
\author{Do Duc Hanh}       
\maketitle

\centerline{\it C/o: Institue of Mathematics, National Centre for Science}\centerline{\it 
 and Technology, P. O. Box 631, Bo Ho, 10.000, Hanoi, VietNam.}
\centerline{e-mail: hanhmath@yahoo.com}
\begin{abstract}In this article we describe the coadjoint orbits of  SL(2,$\mathbb R$). After choosing polarizations for each orbits, we pointed out the corresponding quantum coadjoint orbits and therefore unitary representations of SL(2,$\mathbb R$) via deformation quantization.
\end{abstract}
\section{Introduction}
Let us recall that quantization is a process associating to each Poisson manifold M a Hilbert space H of so-called quantum states, to each classical quantity f$\in \Ci (M)$ a quantum quantity Q(f) $\in \frak L (H)$, i.e., a continuous, perhaps unbounded, normal operator which is auto-adjoint if f is a real-valued function  such that
$$Q(\{f,g\})=\frac{i}{\hbar} [Q(f), Q(g)],$$
$$ Q(1)=Id_H.$$
There are some approachs to this problem, such as Feynman path integral quantization, pseudo differential operator quantization, geometric quantization, etc...In Fedosov deformation quantization, the quantization is considered as the deformation of the structure of the Poisson algebra of classical observables via a family of associated algebras indexed by the so-called deformation parameter rather than a radical change in the nature of the observables.

It is interesting to contruct quantum objects corresponding to the classical ones. It is well-known that the coadjoint orbits are almost all the classificatied flat G-symplectic manifolds. A natural question is to associate to  coadjoint orbits some quantum systems called {\it quantum coadjoint orbits}. Following  Kontsevich'result, every Poisson structure can be quantized. However, this quantizating is only formal and it is difficult to calculate exacly the  corresponding quantum objects and representations in concrete cases. Recently,  Do Ngoc Diep and Nguyen Viet Hai, in \cite{Diep4}, \cite{Diep5}, described the quantum coadjoint orbits and representations of MD and $MD_4$ groups. However, the problem for SL(2,$\mathbb R$) is still open. Although all the irreducible unitary representations of SL(2,$\mathbb R$) are well-known, the correspondence of them with coadjoint orbits are not yet clarified. In this paper, we shall use  Fedosov deformation quantization to find out \s -product formulae and representation of SL(2,$\mathbb R$). 
The algebras of smooth functions on coadjoint orbits of SL(2,$\mathbb R$), deformed by exactly computed \s-products give us series of quantum coadjoint orbits: quantum elliptic hyperboloids, quantum upper (lower) half-hyperboloids, quantum upper (lower) cones, etc...These quantum objects, as we know, appear here for the first time.

The paper is organized as follows. We describe coadjoint orbits in $\S$2. In $\S 3$ we compute for each coadjoint orbit a polarization. The deformation \s-products are computed in $\S$4 and in the last section $\S$5, we show the relation with the unitary dual of SL(2,$\mathbb R$).

For notation, we refer readers to \cite{Kiri1} or \cite{Diep2}, \cite{Diep4}, \cite{Diep5}.
\section{Coadjoint orbits of SL(2,$\mathbb R$)}

Recall that SL(2,$\mathbb R$) is a Lie group with Lie algebra consisting of 2 by 2 matrices with trivial trace. It admits a natural basis of three generators:
\begin{equation*}
H= 
\left(\begin{array}{cc}
1 & 0 \\
0 & -1\end{array}\right),                           
X= 
\left(\begin{array}{cc}
0 & 1 \\
1 & 0 \end{array}\right) ,                          
Y= 
\left(\begin{array}{cc}
0 & 1 \\
-1 & 0 \end{array}\right) ,                           
\end{equation*}
subject to relations: [H,X]=2Y,     [H,Y]=2X,        [X,Y]=-2H. 
Denote by $X^*,H^*,Y^*$ the dual basis of \gs. Because the Killing form is non-degenerate, we can identity $\frak g$ with $\frak g ^*$ in such a way that $\hat X (Y)=\frac{1}{4} B(X,Y)=\frac{Tr(ad X.adY)}{4}$. This isomorphism maps X into $2X^*$, H into $2H^*$, Y into $-2Y^*$.

Naturally, the coadjoint action of SL(2,$\mathbb R$) on  $\g^*$ is given by:
$$\langle K(g) F,Z \rangle =\langle F,Ad(g^{-1})Z \rangle \quad \forall F\in \g^*,g\in G \mbox{, and } Z \in \g.$$
where  $\mathfrak{g}$ is a G-space víi Ad-action. However, there is a natural isomorphism of G-spaces.
\begin{pro}: Operator $X\mapsto \widehat X$ is an smooth G-equivariant isomorphism between G-spaces.
In another words,
$\widehat{Ad(g)X}=K(g)\widehat{X}.$
\end{pro}
It is well-known that GL(2,$\mathbb R$) is a direct product of SL(2,$\mathbb R$) and $R^*=R\backslash(0)$, and therefore each B $\in GL(2,\mathbb R)$ can be decomposed as the product of an element from SL(2,$\mathbb R$) and $\la \left(\begin{array}{cc}1 & 0\\0& 1\end{array}\right)$ or $\la \left(\begin{array}{cc}0 &1 \\1 &0\end{array}\right)\mbox{ with } \la \in R^*_+$.

 Due to the equivariant isomorphism of $\frak g$ with Ad-action and  \gs with K-action, we study the adjoint orbits in place of coadjoint orbits of \gs. It is well-known that every matrix $B\in sl(2,\mathbb R)$ can be reduced to one of the following normal forms:
\begin{equation*}
\left(\begin{array}{cc}
0 & \la\\
-\la& 0\end{array}\right)  ,
\left(\begin{array}{cc}
0 & -\la\\
\la& 0\end{array}\right),  
\left(\begin{array}{cc}
0 & 0 \\
1 & 0\end{array}\right)   ,
\left(\begin{array}{cc}
0 & 1\\
0& 0\end{array}\right)    ,
\left(\begin{array}{cc}
\lambda & 0 \\
0 & -\lambda \end{array}\right)  ,      
\left(\begin{array}{cc}
0 & 0 \\
0 & 0 \end{array}\right),
\end{equation*} 
We obtain the following description of the geometry of coadjoint orbits which is folklore but we could not locate a precise computation from research literature.
\begin{theo}Each coadjoint orbit of SL(2,$\mathbb R$) is one of the forms:
\begin{enumerate}
\item[(a)] Elliptic hyperboloid: $\Omega^1_{\lambda}$=$\{\,  2xX^*+2hH^*-2yY^* \mid x^2+h^2=y^2+\lambda^2, \lambda\not=0\},$\\
\item[(b)] Upper half-cones: $\Omega^2_+$=$\{\,  2xX^*+2hH^*-2yY^* \mid x^2+h^2=y^2, y>0\},$\\
    $\quad$Lower half-cones: $ \Omega^2_-$=$\{\,  2xX^*+2hH^*-2yY^* \mid x^2+h^2=y^2, y<0\},$\\
    $\quad$One point: $\Omega^2_0$=$\{\, 0\},$\\
\item[(c)]Upper half-hyperboloid: $\Omega^3_+$=$\{\,  2xX^*+2hH^*-2yY^* \mid x^2+h^2=y^2-\lambda^2, y>0\},$\\
$\quad$Lower half-hyperboloid: $\Omega^3_-$=$\{\,  2xX^*+2hH^*-2yY^* \mid x^2+h^2=y^2+\lambda^2, y<0\}.$\\
\end{enumerate}
\end{theo}
{\bf Proof}.
We describe the geometry of adjoint orbits corresponding to  $\Omega^1_{\lambda}$, $\Omega^2_-$ and $ \Omega^3_{\la,+}$. The case of other orbits can be analogously treated.
The adjoint orbit corresponding to $\Omega^1_{\lambda}$ contains 
$ \left(\begin{array}{cc}
\lambda & 0 \\
0 & -\lambda \end{array}\right)$. By a direct computation, for $S=
\left(\begin{array}{cc}
u & v \\
s & t \end{array}\right)\in SL(2,\mathbb R)$, we have 
\begin{equation*}
\left(\begin{array}{cc}
h & x+y \\
x-y & -h \end{array}\right)=
S \left(\begin{array}{cc}
\lambda & 0 \\
0 & -\lambda \end{array}\right)S^{-1}=\left(\begin{array}{cc}
\lambda(ut+sv) & -2\la uv \\
2\la st & -\la(ut+sv) \end{array}\right).
\end{equation*}
Hence, $\frac{h}{\la}=ut+sv$, $\frac{x+y}{\la}=-2uv$, $\frac{x-y}{\la}=2st$ and therefore, $\frac{x^2-y^2}{\la^2}+\frac{h^2}{\la^2}=-4uvst+(ut+sv)^2=(ut-sv)^2=1.$
Moreover, the coadjoint orbit  containing $2\la H^*$ is \\
 $\{2xX^*+2hH^*-2yY^*\mid x^2+h^2-y^2=\la^2\}$.\\ It is exactly the elliptic hyperboloid.
The adjoint orbit corresponding to $\Omega^3_{\la,-}$ containing 
$ \left(\begin{array}{cc}0 & -\la \\ \la & 0 \end{array}\right) $.  
By a direct computation, for $S=\left(\begin{array}{cc}u & v \\s & t \end{array}\right)\in SL(2,\mathbb R)$, we have 
\begin{equation*}\left(\begin{array}{cc}
h & x+y \\
x-y & -h \end{array}\right)
=S \left(\begin{array}{cc}
0&-\la  \\
\la & 0 \end{array}\right)S^{-1}=\left(\begin{array}{cc}
\lambda(vt-us) & -\la (u^2+v^2) \\
\la (s^2+t^2) & -\la(vt-us) \end{array}\right).
\end{equation*}
Hence, $\frac{\ds h}{\ds \la}=vt+us$, $\frac{\ds x+y}{\ds \la}=-(u^2+v^2)$, $\frac{\ds x-y}{\ds \la}=s^2+t^2$.
And therefore, $\frac{\ds x^2-y^2}{\ds \la ^2}+\frac{\ds h^2}{\ds \la^2}=1 $ for $0\geq x+y$, $x-y\geq 0$.
Moreover, the coadjoint orbit  containing $2\la Y$ is 
$\{2xX^*+2hH^*-2yY^*\mid x^2+h^2=y^2-\la^2,y<0\}$
It is exactly one of the two connected components of a the elliptic  hyperboloid. 
Let us consider the adjoint orbit corresponding to $\Omega^2_-$ containing $ \left(\begin{array}{cc}
0 & 0 \\
1 & 0 \end{array}\right)$\\
By direct computation, for $S \in SL(2,\mathbb R)$, we have:
\begin{equation*}
\left(\begin{array}{cc}
h & x+y \\
x-y & -h \end{array}\right)=
S \left(\begin{array}{cc}
0 & 0 \\
1 & 0 \end{array}\right)S^{-1}\left(\begin{array}{cc}
vt & -v^2 \\
t^2 & -vt \end{array}\right).
\end{equation*}
Hence, h=vt, $x+y=-v^2$, $x-y=t^2$.
And therefore $x^2+h^2-y^2=0, 0\geq x+y, x-y \geq 0$. Note that $(x, h, y)\not=(0, 0, 0)$. The coadjoint orbit  containing $X^*+Y^*$ is\\
$  \{2xX^*+2 hH^*-2 yY^*\mid x^2+h^2=y^2, y>0\} $
It is really the upper half-cones.\\
The theorem is proved.
\section{Complex Polarirations of K-orbits of SL(2,$\mathbb R$)}                                                                                                                                                                                                                                                                                                                                                                                                                                                                                                                                                                                                                                     
Before quantizing coadjoint orbits we do first describe some polarizations on orbits. Let us recall some basis concepts concerning polarization, see \cite{Diep2}.

Let G be a Lie group. A complex polarization of orbit $\Omega _F$ at  $ F \in \Omega _F$ is a quadriple of  ($\eta, \frak h , U, \rho$) such that:
\begin{enumerate}
\item $\eta$ is a subalgebra of the complex Lie algebra  $\g_C=g\underset{R}{\otimes} C$ containing $\g_F$.
\item The subalgebra $\eta$ is invariant under the action of all the operators of type $Ad_{\g_C} x$, $x \in G_F.$
\item The vector space $\eta + \bar \eta$  is complexification of real subalgebra Lie $\frak m=(\eta+\bar \eta)\cap \g.$ \item All subgroup $M_0,H_0,M,H$ are closed, where, by definition $M_0$ (resp., $H_0$) is the connected subgroup of G with Lie algebra m (resp., $\h :=\eta \cap \g$) and M:=$G_F.M_0$, H:=$G_F.H_0$.
\item  U is an irreducible representation of $H_0$ in some Hilbert space H such that: 1. The restriction $U\mid _{G_F\cap H_0}$ is some multiple of $\chi_F$ 
where by definition $\chi_F(exp X)\mid_{(G_F)_0\cap H_0}:=exp(2\pi \sqrt{-1} \langle F,X \rangle)$; 2. The Nelson condition is satisfied. See \cite{Diep2},  $10.5$, theorem 3.
\item The Pukanszky condition is satisfied: $F+\eta^\bot \subset \Omega_F$, see \cite{Kiri1}, $\S 15.3$
\end{enumerate}
Denote by $\rho$ the one dimension reprentation $2\pi \sqrt{-1} \langle F,X \rangle$ of Lie algebra $\eta$.\\Let $\Ci (G, \eta, H, \rho, U)$ be the set of common solutions of 
$$f(hg)=U(h). f(g)$$
$$(L_X -\rho (X))f=0\quad X\in \eta $$
\begin{re}
The condition 5 and 6 are often included in order to obtain irreducible representations.
\end{re}
In this section, we establish complex polarization for K-orbits.
\subsection{Polarization of $\Omega ^1 _\lambda$}
Let us consider a point $\hat F=2\la H^*\in \Omega^1_\la $, the complex subalgebra $\eta=\langle H,X+Y\rangle_C$. The representation  $U=e^{2\pi i \langle F,.\rangle}$ of $\h=\eta\cap \g$ can be extended to $H=H^0\cup \varepsilon H^0$ as $U(\varepsilon )=\pm 1$. Let $\rho $ be the 
natural extension of $dU$ to $\eta$
\begin{pro}
$(\eta, \rho, U)$ is a polarization of $\Omega^1_\la$.
\end{pro}
{\bf Proof}. It is easy to see that the stabilizer $G_F =\left\{ \left( \begin{array}{cc}
                       a& 0 \\
                          0 & a^{-1}
                        \end{array} \right) \right\}$
 consists of two connected components corresponding to  $a>0$ and $a<0$. Obviously, its Lie algebra is $\g _F=\langle H\rangle $. The  Ad-orbit passing through $F=\la H$  contains two lines $\{ F+t(X\mp Y)\}$. Clearly, these lines are the images of ones $\{\hat F+t(X^*\pm Y^*)\}$ passing through $\hat F$ on $\Omega^1_\la$ under the isomorphism generated by Killing form. 
Chose $\eta=\langle H, X+Y\rangle _C$. We can see Pukansky satisfied.
Note that [H, X+Y]=2(X+Y) so $\eta $ is a invariant Lie algebra under Ad-action of  $G_F$. We also deduce $\frak h=\eta \cap \g=\frak m=\langle H, X+Y\rangle , \bar \eta=\eta, \frak m_C=\eta+\bar \eta=\eta$. Chose $\rho(A)=2\pi i \langle \hat F, A\rangle $ with $A\in \eta$ is holomorphic representation of $\eta$. We have, $\rho (aH+b(X+Y))=4\pi i \la  a$. Because $G_F$ has two connected components,
$H=G_F . H^0 =\left\{ \left( \begin{array}{cc}
                       \alpha & \beta \\
                          0 & \alpha^{-1}
                        \end{array} \right)\mid \alpha \not= 0 \right\}$.\\
By an exact computation, we have\\
$exp \left(\begin{array}{cc}
a & b \\
0 &-a \end{array}\right)=exp\Biggl(a.     \left(\begin{array}{cc}
1 & 0 \\
0 & -1 \end{array}\right)+b\left(\begin{array}{cc}
0 & 1 \\
0 & 0 \end{array}\right)\Biggr)=\left(\begin{array}{cc}
e^a & b(\frac{e^a-e^{-a}}{2}) \\
0 & e^{-a} \end{array}\right)    $.\\
Thus,
$U\left(exp
 \left(\begin{array}{cc}
a & b \\
0 &-a \end{array}\right)\right)=e^{4\pi i \la a}$
or $ U\left( \begin{array}{cc}
                       \alpha & \beta \\
                          0 & \alpha^{-1}
                        \end{array} \right)=\alpha^{4 \pi i \la}$ for all $\la>0$.
On the other hand, $ H=H^0\cup \left( \begin{array}{cc}
                      -1 & 0 \\
                        0 &-1
                        \end{array} \right). H^0$,
and so we can extend U onto H following 
$U\left(\begin{array}{cc}
-1 & 0 \\
0 & -1 \end{array}\right)=\pm I$. Corresponding to characters of $H/ H^0=\mathbb Z_2$, we obtain thus two unitary representations of  H: 
$U \left(\begin{array}{cc}
\alpha & \beta \\
0 & \alpha^{-1} \end{array}\right)=|\alpha|^{4 \pi i\la}$
 and $U\left(\begin{array}{cc}
\alpha & \beta \\
0 & \alpha^{-1} \end{array}\right)=|\alpha|^{4 \pi i\la}. sgn(\alpha)$.
\subsection{Polarization of orbit $\Omega^2_+$ }
                                 
Let us consider a point $\hat F=X^*-Y^*\in \Omega^2_+ $, the complex subalgebra $\eta=\langle H, X+Y\rangle_C$. The representation  $U=e^{2\pi i \langle F,.\rangle}$ can be extended to $H=H^0\cup \varepsilon H^0$ as $U(\varepsilon )=\pm 1$. Let $\rho $ be the 
natural extension of $dU$ to $\eta$.
\begin{pro}
$(\eta, \rho, U, \rho)$ is a polarization of $\Omega^2_+$.
\end{pro}
{\bf Proof}. It is easy to see that the stabilizer $G_F =\left\{ \left( \begin{array}{cc}a& b \\  0 & a                        \end{array} \right) \right\};a \in\{-1, 1\}$ consists of two connected components corresponding to $a>0$ and $ a<0$ with Lie subalgebra $\g _F=\langle X+Y\rangle $. 
Chose  $\eta=\langle H, X+Y\rangle _C$. 
Due to [H, X+Y]=2(X+Y), $\eta $ is a invariant Lie algebra under the  Ad action of  $G_F$. As known
$\eta^\bot=\langle X^*-Y^*\rangle $: functional on $\frak g$ such that vanishes on $\eta$ when extented to complexification of \g.
We also imply $\frak h=\eta \cap \g=\langle H, X+Y\rangle , \bar \eta=\eta=\frak m_C$ and 0 is the one-dimension representation of $\eta$. 
Naturally,
$ H=H^0\cup \left( \begin{array}{cc}
                      -1 & 0 \\
                        0 &-1
                        \end{array} \right). H^0$ 
and $\left(\begin{array}{cc}
-1 & 0 \\
0 & -1 \end{array}\right)^2=I$.
It follows $U\left(\begin{array}{cc}
-1 & 0 \\
0 & -1 \end{array}\right)=\pm I$.
Following to characters of $H/ H^0$ we obtain two unitary representations of  H:\\
$U \left(\begin{array}{cc}
\alpha & \beta \\
0 & \alpha^{-1} \end{array}\right)=1
$ and $ U\left(\begin{array}{cc}
\alpha & \beta \\
0 & \alpha^{-1} \end{array}\right)=\mbox{sgn}(\alpha)$.
By analogy, we obtain the same result for $\Omega^2_-$.
\subsection{Polarization for $\Omega^3_{\lambda, +}$ }
Let us consider a point $\hat F=2 H^*\in \Omega^3_{\lambda, +}$, the complex subalgebra $\eta =\langle Y,X+iH \rangle_C$. Because of the fact that the stabilizer SO(2,$\mathbb R$) of $\hat F$ is not simply connected, $U=e^{2\pi i \langle F,.\rangle}$ can be extented to H only if the orbit is integral.
\begin{pro}
$(\eta, \rho,U,\rho)$ is a polarization of $\Omega^3_{\la,+}$ and this orbit is integral if and only if $\la$ is of the form $\la=\frac{k}{8}$.
\end{pro}{\bf Proof}. It is trivial that the stabilizer
$G_F =SO(2, R)$ with Lie algera $\g _F=\langle Y\rangle $ is connected but not simply connected. By choosing 
 $\eta=\langle Y, X+iH\rangle _C,\frak m_C=\g, \frak h=\eta \cap \g $
, $\eta$ admits an one-dimension representation
$\rho \left(\begin{array}{cc}
-ia&a+b\\
-a+b&ia
\end{array}\right)=-4\pi i \la a$,
which has the restriction on $\frak h$,   
$\rho \left(\begin{array}{cc}
0&a\\
-a&0
\end{array}\right)=-4\pi i \la a$. 
On the other hand, \\
\begin{equation*}
exp \left(\begin{array}{cc}
0&a\\
-a&0
\end{array}\right)=\left(\begin{array}{cc}
\cos  a &\sin a\\
-\sin a &\cos a
\end{array}\right).
\end{equation*}
Thus
\begin{equation*}
U\left(\begin{array}{cc}
\cos  a&\sin a\\
-\sin  a&\cos a
\end{array}\right)=e^{-4\pi i \la a}.
\end{equation*}
Because SO(2,$\mathbb R$) is not simply connected, U may not exist.
The nesesary and sufficient condition is $\la=\frac{k}{8}$.
The orbit $\Omega^3_{\la, -}$ can be treated analogously and we gain the same result.
A corollary of polarization for all co-adjont orbits is the representation of SL(2,$\mathbb R$) on the Hilbert space of partial holomophic, square- integrable sections of induced vector bundle. See e.g \cite{Lang}, \cite{Diep2}.  We follow another approach by deformation quantization.
\section{Quantum coadjoint orbits of SL(2,$\mathbb R$)}
We shall work from now on for the fixed coadjoint orbit $\Omega^1_{\lambda}$. Following the scheme from \cite{Diep4},\cite{Diep5}, first we study he geometry of this orbit and introduce some canonial coordinates in it. It's well known that coadjoint orbits are isomorphism to the homogeneous spaces $G/G_F$ which are symplectic manifolds.  We'll introduce a coordinate system on this orbit and it turns out to be a Darboux one. Each $A \in \g$ can be considered as a linear functional $\A$ on coadjoint orbits, as a subset of $\g^*, \A (F)=\langle F,A \rangle$. It is also well known that this function is just the Hamilton function associated with the Hamiltonian vector field $\xi_A$ generated by A following the formula:
$$\xi_A (f) (x)=\frac {d}{d t}f( x \exp{(tA))}  \mid _{t=0}$$
The Kirillov form $ \omega_F$ is defined by the formula
$$\omega_F (\xi _A,\xi _B )=\langle F,[A,B]\rangle$$
It is known as the flatness of the coadjoint orbits that the correspondence $A \mapsto \A$ is a Lie homomorphism.
Motivated by the contructed polarizations, $\Omega^1 _\la$ can be parameterized as
\begin{equation*}
\begin{cases}
x=M(p, q)=p\cos (q)-\lambda \sin (q);\\
h=N(p, q)=p\sin (q)+\lambda \cos (q);\\
y=P(p, q)=p;
\end{cases}
\end{equation*}
M, N, P satisfy
\begin{equation}
M_q=-N;N_q=M;M_p=\cos (q);N_p=\sin (q);M. \cos (q)+N. \sin (q)=p;\label{E36}
\end{equation}
Let us consider the mapping $\psi:(p,q)\mapsto 2M(p,q)X^*+2N(p,q)H^*-2P(p,q)Y^*$
Clearly, $(R^2, \Omega ^1 _\lambda, \psi)$
is an universal covering space.
\begin{pro}
$\psi$ is a symplectomophism and Hamiltonian \A in coordinates (p, q) is of the form: 
$$\A(F)=\langle F, A\rangle =(2 a_1 \cos  q+2b_1 \sin q-2 c_1)p+(-2a_1 \sin q+2 b_1 \cos q)\lambda$$
\end{pro}
{\bf Proof}:$\quad$
Each F$\in \Omega^1_\la $ is of the form $2MX^*+2NH^*-2PY^*$. From this it folllows that the Hamiltonian function generated by invariant vector field $\xi_A$ is
$$\A (F)=\langle F, A\rangle  =2 a_1 M+2 b_1 N-2 c_1 P.$$
It follows therefore
$$\A(F)=2a_1(p\cos q-\lambda \sin q)+2b_1(p \sin q+\lambda \cos  q)-2c_1p.$$
On $R^2$ there are two symplectic structures: the first one is the Kirillov form induced by mapping $\psi$ and the second is the canonical symplectic form $dp\wedge dq$. We prove their coincidence by observing their values at invariant vector fields are equal.\\
 Note that $\omega_F(\xi_A, \xi_B)=\langle F, [A, B]\rangle  \\= \langle 2MX^*+2NH^*-2PY^*, 2(b_1c_2-b_2c_1)X+2(c_1a_2-c_2a_1)H-2(a_1b_2-a_2b_1)Y\rangle $
$=4M(b_1c_2-b_2c_1)+4N(c_1a_2-c_2a_1)+4P(a_1b_2-a_2b_1).$\\
On the other hand,\\
$(dp\wedge dq) (\xi_A,\xi_B)=\{\A, \B\}=\frac{\p \A}{\p p}\frac{\p \B}{\p q}-\frac{\p \A}{\p q}\frac{\p \B}{\p p}\\=4(b_1c_2-b_2c_1)N_q+4(c_1a_2-c_2a_1)(-M_q)+4(a_1b_2-a_2b_1)(M_pN_q-N_pM_q).$\\
Then $\omega_F(\xi_A, \xi_B)=(dp\wedge dq) (\xi_A,\xi_B)$.\\
The theorem is therefore proven.\\
\begin{re}
The case of diffenrent orbits can be treated similarly with a small correction.
With the orbits $\Omega^3_{\la,+}$ and $\Omega^3_{\la,+}$, clearly we can't find out a affine subspace of a half dimensions, thus there can't exist a coordinate as above. However, a good approach is considering the complexification of orbits and we  obtain $(C\times C,\Omega_\lambda^1,\psi)$ as universal complex symplectic covering space, only by replacing $\la$ by $i\la$. The orbits $\Omega^2_+$, $\Omega^2_+$ can be viewed as a part of the case $\Omega^1_{\la,+}$ and $\Omega^3_{\la,+}$ when $\la =0$. 
\end{re}
>From now, because of the similarity, we'll deal mainly with the orbits
$\Omega_\lambda^1$. The other orbit can be treated with a simple modification.
\begin{theo}
With A, B $\in \g$, the Moyal \s-product satisfies $$i\A \s i \B-i\B\s i\A=i\widetilde{[A, B]}$$
\end{theo}
{\bf Proof}:\hspace*{20pt}
Consider two arbitrary elements $A=a_1 X+b_1 H+c_1 Y, B=a_2 X+b_2 H+c_2 Y \in \g$ , 
By the Moyal-Weyl formular,
$$i\A \s i\B=\sum_{k=0}^{\infty}P^k(i\A, i\B). \frac{1}{k!}(\frac{1}{2i})^k,$$
with $P^k(i\A, i\B)=-\wedge^{i_1j_1}\wedge^{i_2j_2}\cdots \wedge^{i_kj_k} \p_{i_1i_2\cdots i_k}\A \p_{j_1j_2 \cdots j_k} \B$\\
It's easy, then, to see that:
$$P^0(i\A, i\B)=-\A. \B,$$
$$P^1(i\A, i\B)=-(\wedge^{12}\frac{\p \A}{\p p}. \frac{\p \B}{\p q}+\wedge^{21}\frac{\p \A}{\p q}. \frac{\p \B}{\p p})=-\{\A, \B\}, $$
By proposition 4.1, \A, \B are linear functions of p. Thus for $k \geq 2$, we have
\begin{align*}
P^2(i \A, i\B)&=-(\wedge^{12}\wedge^{12}\A_{pp}\B_{qq}+\wedge^{21}\wedge^{21}\A_{qq}\B_{pp}&\\&+
\wedge^{12}\wedge^{21}\A_{pq}\B_{qp}+\wedge^{21}\wedge^{12}\A_{qp}\B_{pq}=-2\A_{pq}\B_{qp}.&
\end{align*}
 $P^2(i \A, i\B)=P^2(i \B, i\A)$,.
Therefore\\
$P^k(i\A, i\B)=-\wedge^{i_1j_1}\wedge^{i_2j_2}\cdots \wedge^{i_kj_k} \p_{i_1i_2\cdots i_k}\A \p_{j_1j_2 \cdots j_k} \B =0\hspace*{20pt} \forall k \geq 3.$\\
We get\\
$i\A \s i \B-i\B\s i\A=(P^1(i\A, i\B)-P^1(i\B, i\A))\frac{1}{2i}+(P^2(i\A, i\B)-P^2(i\B, i\A))(\frac{1}{2i})^2. \frac{1}{2!}=i\{\A, \B\}=i\widetilde {[A, B]}$. \\
The theorem can be proved analogously on $\Omega^2_+$, $\Omega^2_-$ and $\Omega^3_{\lambda, C} .$ 
\begin{re}
Consider the canonial representation of quantum algebra $(\Ci(\Omega),\s)$ on itself which is a FrÐchet Poisson algebra  by left \s-multiplication defined by:
$$l_f:\Ci(\Omega) \rightarrow \Ci(\Omega),$$
$$ g \mapsto f\s g.$$
Then,  $\Ci(\Omega)$ can be viewed as a algebra of pseudo-diiffefential operators on $\Ci (\Omega)$.    
On the other hand, the corespondence $ A \mapsto \A$  is a Lie algebra homomorphism. Thus, we can consider the repersentation of Lie algebra sl(2,$\mathbb R$) on dense subspace $L^2(\mathbb R\times [0,2\pi), \frac {dpdq}{2 \pi} )^\infty $ of smooth functions by left \s-multiplication  by $i \A\s$. This representation is then extended to the whole space $L^2(\mathbb R\times SO(2,\mathbb R),\frac{dp.dq}{2 \pi})$ by \cite{Arnal1}. We study now the convergence of the formal power series. In order to do this, we look at the \s-product of i\A as the \s -product of symbols and define the differential operators corresponding to i\A . It is easy to see that the resulting correspondence is a representation of $\frak g$ by pseudo-differential operators. 
\end{re}
On $\Omega ^\lambda_1=\{2xX^\ast+2hH^\ast-2yY^\ast \mid x^2 +h^2=y^2+\lambda ^2\}$ the following results hold:
\begin{lem}
\item{(1)}$ \fp(\p _p \fp ^{-1} (f))=i^{-1} (x.f),$
\item{(2)}$\fp (p.\fp^{-1} f)=i \p_x (f)$,
\item{(3)}$P^k(\A, \fp^{-1} (f))=k(-1)^{k-1}\A_{q\cdots q p}\p_{p \cdots pq}\fp^{-1}(f)$
$+(-1)^k \A_{q\cdots q}\p_{p \cdots}\fp^{-1}(f)$.
\end{lem}
{\bf Proof}. The two first formulas are well-known from the theory of Fourier transforms.
If  $k \geq 2$ then by theorem 4.1 it implies \A is a linear function of p. Because one of the coordinates is linear, if two of index $i_1, i, _2, \cdots, i_k$ equals to 1 then $\p_{i_1, i, _2, \cdots, i_k} \A=0$. Therefore, for all $k \geq 2$:
\begin{multline*}
P^k(\A, \fp^{-1} (f))=\wedge ^{i_1, j_1}\wedge ^{i_2j_2}\cdots \wedge^{i_kj_k}\A_{i_1\cdots i_n}\p_{j_1\cdots j_n}\fp^{-1}(f)\\
=\sum \wedge^{21}\cdots\wedge^{12}\cdots \wedge^{21}\A_{q\cdots p\cdots q}\p_{p\cdots q\cdots p}\fp^{-1}(f)+\v^{21}\cdots \v^{21}\A_{\qq}\p_{\pp}\fp^{-1}(f).
\end{multline*}
It is clear that $\v^{-1}=
\left(\begin{array}{cc}
1 & 0 \\
0 & -1\end{array}\right)  ,$ 
So we get $\v^{12}=1, \v^{21}=-1$ . It deduces\\
$P^k(\A, \fp^{-1} (f))=k(-1)^{k-1}\A_\qpq\p_\pqp\fp^{-1}(f)+(-1)^{k-1}\A_\qq\p_\pp\fp^{-1}(f).$\\
With k=0 hay k=1, clearly, the lemma is also satified. 
Apply this lemma, we have the followimg theorem: 
\begin{theo}If we set s=$q-\frac{x}{2}$, for each compactly supported smooth function f$ \in \Ci_c(\mathbb R^2)$  we have\\
$\hat l_A=\fp\circ l_A\circ \fp^{-1}=(a_1 \cos s+b_1\sin s-c_1)\p_s+(-a_1\sin s+b_1\cos s)(2\lambda i+1)$
\end{theo}
{\bf Proof}.\hspace{20pt}\\
Following the Moyal-Weyl formula,  we have
$$\hat l_A(f)=\fp\circ l_A\circ \fp^{-1}(f)=
=i\fp(\sum_{k=0}^\infty  (\frac{1}{2i})^k. 
\frac{1}{k!}P^k(\A, \fp^{-1}(f)),$$
By the lemma, 
\begin{align}\notag
\hat l_A(f)&=i\fp \Bigl (\sum_{k=0}^\infty (\frac{1}{2i})^k\frac{1}{k!}. (-1)^{k-1}. k. \A_\qpq\p_\pqp\fp^{-1}(f)&\\\notag 
&+ \sum_{k=0}^\infty (\frac{1}{2i})^k\frac{1}{k!}. (-1)^k. \A_\qq\p_\pp\fp^{-1}(f)\Bigr )=I+J,&\notag
\end{align}
Note the fact that \A is a linear function of p. Therefore $\A_\qpq $ is a function of only variable p. 
\begin{align}
I&=i\fp(\sum_{k=0}^\infty (\frac{1}{2i})^k\frac{1}{k!}. (-1)^{k-1}. k. \A_\qpq\p_\pqp\fp^{-1}(f))&\\\notag
&=i\sum_{k=0}^\infty (\frac{1}{2i})^k\frac{1}{k!}. (-1)^{k-1}. k. \fp(\A_\qpq\p_\pqp\fp^{-1}(f))&\\\notag
&=i\sum_{i=1}^\infty(\frac{1}{2i})^k\frac{1}{(k-1)!}. (-1)^{k-1}. (ix)^{k-1}\A_\qpq\p_p f&\\\notag
&=\frac{1}{2}\p_p\A(q-\frac{x}{2})\p_q(f).& \notag
\end{align}
Set \A=p. M+N, where M, N depend only q, by exact computations, we have
\begin{align*}
J&=i\sum _{k=0}^\infty (\frac{1}{2i})^k\frac{1}{k!}. (-1)^k. \fp ((p. M_\qq +N_\qq). \fp^{-1}(f))&\\\notag
&=i\sum_{k=0}^\infty( \frac{i}{2})^k \frac{1}{k!}((i. \p _x M^{(k)}+N^{(k)}). (ix)^k. f)&\\\notag
&=i\sum_{k=0}^\infty (\frac{i}{2})^k \frac{1}{k!}. i\p_x. M^{(k)}(q). (ix)^k. f+i\sum_{k=0}^\infty (\frac{i}{2})^k\frac{1}{k!}N^{(k)}(q)(ix)^k. f&\\\notag
&=-\sum_{k=0}^\infty(-\frac{x}{2})^k. \frac {M^{(k)}}{k!}\p_xf-\sum_{k=0}^\infty(-\frac{1}{2})k. \frac {M^{(k)}}{k!}. k. x^{k-1}\p_xf+i\sum_{k=0}^\infty(-\frac{x}{2})^k. \frac {N^{(k)}}{k!}\p_xf&\\\notag
&=\frac{1}{2}\sum _{k=0} ^\infty (-\frac{x}{2})^k. \frac{M^{(k+1)}(q)}{k!}. f-M(q-\frac{x}{2}). \p _x f+iN(q-\frac{x}{2}).  f&\\\notag
&=\frac{1}{2}. M'(q-\frac{x}{2}). f-M(q-\frac{x}{2}). \p _x f+iN(q-\frac{x}{2}).  f.&\notag
\end{align*}
Finally, we have the exact formular of corresponding quantized operator:
\begin{align}\notag
\hat l_A(f)&=\frac{1}{2}\p_p\A(q-\frac{x}{2})\p_q(f)+\frac{1}{2}. M'(q-\frac{x}{2}). f+-M(q-\frac{x}{2}). \p _x f+iN(q-\frac{x}{2}).  f&\\\notag
&=M(q-\frac{x}{2})(\frac{1}{2} \p _q-\p_x)f+\frac{1}{2}. M'(q-\frac{x}{2}). f+iN(q-\frac{x}{2}).  f.&\\\notag
\end{align}
Put $q-\frac{x}{2}=s ;q+\frac{x}{2}=t$, it follows $\p _s=\p _q-2\p _x$.
Recall that\\ $\A(F)=2a_1(p\cos q-\lambda \sin q)+2b_1(p \sin q+\lambda \cos  q)-2c_1p$.
$M(q)=2a_1(p\cos q-\lambda \sin  q)+2 b_1(p\sin  q+\lambda \cos  q)-2 c_1,$\\
$N(q)=-2\lambda a_1 \sin q+2\lambda b_1 \cos  q,$\\
$M'(q)=\frac{N(q)}{\lambda}.$\\
Therefore,\\
$\hat l_A(f)=\frac{1}{2}M(s). \p _s f+(\frac{N(s)}{2. \lambda}+iN{s})f$
$=(a_1 \cos  s+b_1 \sin  s-c_1)\p _s+(-a_1 \sin  s+b_1 \cos  s)(2\lambda i+1).$\\
The theorem is proved.\\
By analogy, we get the same results for all two dimesion coadjoint orbits.
Note that, following the virtual of the polarizations chosen for orbits, we obtain the representation of sl(2,$\mathbb R$) on $L^2 $-space on SO(2,$\mathbb R$). 
\section{Relation with unitary dual of SL(2,$\mathbb R$)}
We recall some basic results of contructing unitary dual of SL$(2,\mathbb R)$ by the classical methods, see e.g. \cite{Lang}.\\
Consider  the subgroup $H=\left(\begin{array}{cc}a&b \\0&a^{-1}\end{array}\right)$ associated with one-dimension representation $\rho_s\left(\begin{array}{cc}a&b \\0&a^{-1}\end{array}\right)=a^{s+1}$. Let $\phi_s $ be the induced representation of $\rho_s$ on to SL(2,$\mathbb R$). Clearly, the space of induced vector bundle is isomorphic to the space $H_s$ of function on G satisfies $f(hg)=\rho_s(h).f(g)$ with restriction on K lying on $L^2(K)$, also isomorphic to $L^2(K)$ where $K=SO(2,\mathbb R)\simeq G/H$  

Let T be the representation of G on $\Ci (G)$ defined by $T(g_1)f(g)=f(gg_1)$. The infinitesimal representation of T determined by $L(A)f(g_0 )=\frac{\p}{\p t} T(e^{tA}) f(g_0)\mid _{t=0}$.
By the Iwasawa decomposition, each g of SL(2,$\mathbb R$) can be viewed as the product
$g=\left(\begin{array}{cc}
\frac{ 1}{ \sqrt y} & 0 \\0 & \frac{1}{ \sqrt y}\end{array}\right).     
\left(\begin{array}{cc}
y & x \\
0 & 1\end{array}\right). 
\left(\begin{array}{cc}
\cos (\theta) & \sin (\theta) \\
-\sin (\theta) & \cos (\theta)\end{array}\right). $   \\ 
So a function on G can be viewed as function of x, y, $\theta$. We obtain the explicited fomulars of L as:\\
$L_X=(s+1)\sin 2\theta-\cos  2 \theta \p _\theta,$\\
$L_H=(s+1) \cos  2\theta+\sin  2 \theta \p _\theta,$\\
$L_Y=\frac{\p}{\p \theta},$\\
 From this, by considering the algebraic vector subspaces of $L^2(K)$, it can imply all the irreducible unitary representations of $SL(2,\mathbb R)$ of discrete series, principal series, the complementary series as in \cite{Lang}. In order to prove the equivalence of two approachs, it is enough to show that the corresponding infinitesimal representations of Lie algebra sl(2,$\mathbb R$) are the same.
\begin{theo}The  representations $\hat l$ obtained from deformation quantization are coincided with the infinitesimal representation L of Lie algebra corresponding to discrete series, principal series, the complementary series of SL(2,$\mathbb R$). 
\end{theo}
{\bf Proof}: We know that $f(x, y, \theta)=y^{\frac {s+1}{2}}. f(\theta)$.\\
So, $\frac{\p f}{\p y}=\frac{s+1}{2y}. f$ ,  $\frac{\p f}{\p \theta}=0.$\\
Thus $2y \p _y=s+1$, $\p _x=0$. We obtain the explicited formual of representation: for $A=a_1 X+b_1H+c_1 Y$\\
$L_A=(-a_1 \cos (2 \theta)+b_1 \sin  2\theta+c_1)\p _\theta+(s+1)(a_1 \sin  2 \theta+b_1\cos  2 \theta)$\\
setting $s=2\la$ vµ $-2\theta=s$. Then\\
$L_A=(a_1 \cos  s+b_1 \sin  s-c_1)\p _s+(-a_1 \sin  s+b_1 \cos  s)(2\la+1)=\hat l_A$\\
The proof is therefore achieved.\\
\begin{re}
We demonstrated how irreducible unitary representations of SL(2,$\mathbb R$) could be obtained from deformation quantization. It is reasonable to refer to the algebras of functions on coadjoint orbits with corresponding \s -product as  a quantum ones, namely {\bf quantum elliptic hyperboloids} $(\Ci (\Omega^1_\la),\s_\hbar)$, {\bf quantum elliptic cones} $(\Ci (\Omega^2_\pm),\s_\hbar)$, {\bf two folds quantum hyperboloids } $(\Ci(\Omega^3_\la),\s_\hbar)$ etc.
\end{re}

\end{document}